\documentclass[11pt]{article}
\usepackage{amstext,amssymb,amsmath,amsbsy}

\usepackage{hyperref}
\usepackage{amscd}
\usepackage{amsfonts}
\usepackage{indentfirst}
\usepackage{verbatim}
\usepackage{amsmath}
\usepackage{amsthm}
\usepackage{enumerate}
\usepackage{graphicx}
\usepackage{subfig}
\usepackage{color}
\usepackage[OT1]{fontenc}
\usepackage[latin1]{inputenc}
\usepackage[english]{babel}
\usepackage{amssymb}
\newtheorem{theorem}{Theorem}
\newtheorem{lemma}{Lemma}

\numberwithin{equation}{section}

\newcommand{\proofend}{\hfill $\Box$ }

\newcommand{\supp}{\operatorname{supp}}

\newcommand{\dive}{\operatorname{div}}

\newcommand{\eps}{\varepsilon}

\newcommand{\epss}{s}

\newcommand{\mC}{\mathbb{C}}

\newcommand{\mR}{\mathbb{R}}

\newcommand{\hu}{\hat u}

\newcommand{\bU}{{U}}
\newcommand{\bV}{{V}}

\title{Cloaking using complementary media in the quasistatic regime}

\author{Hoai-Minh Nguyen \footnote{EPFL SB MATHAA CAMA, Station 8,  CH-1015 Lausanne, hoai-minh.nguyen@epfl.ch}  \footnote{The research is partially supported by NSF grant DMS-1201370 and by the Alfred P. Sloan Foundation.}}

\begin{document}

\maketitle 

\begin{abstract} Cloaking using complementary media was suggested by Lai et al. in \cite{LaiChenZhangChanComplementary}. The study of this problem faces two difficulties.  Firstly, this problem is unstable since the equations describing the phenomenon  have sign changing coefficients, hence the ellipticity is lost. Secondly, the localized resonance, i.e., the field explodes in some regions and remains bounded in some others as the loss goes to 0, might appear.  In this paper, we give a proof of cloaking using complementary media  for a class of  schemes inspired from \cite{LaiChenZhangChanComplementary} in the quasistatic regime. To handle the localized
resonance, we introduce the technique of removing localized singularity and apply a  three spheres inequality.  The proof also uses the reflecting technique in \cite{Ng-Complementary}. To our knowledge, this work presents the first proof on cloaking using complementary media. 
\end{abstract}



\section{Introduction}
Negative index materials (NIMs) were first investigated theoretically by Veselago in \cite{Veselago} and were innovated by Nicorovici et al. in \cite{NicoroviciMcPhedranMilton94}  and Pendry in~\cite{PendryNegative}. The  existence of such materials was confirmed by Shelby, Smith, and Schultz in \cite{ShelbySmithSchultz}. The study of NIMs has attracted a lot the attention of the scientific community thanks to their many possible applications. One of the appealing ones is cloaking. There are at least three ways to do cloaking using NIMs. The first one  is based on the concept of anomalous localized resonance discovered by Milton and Nicorovici in \cite{MiltonNicorovici}. The second one is based on plasmonic structures introduced by Alu and Engheta in \cite{AluEngheta}. The last one makes use of  the concept of complementary media and was suggested by Lai et al. in \cite{LaiChenZhangChanComplementary}. 
In this paper, we concentrate on the last method.  

\medskip
The study of cloaking using complementary media faces two difficulties.  Firstly, this problem is unstable since the equations describing the phenomenon  have sign changing coefficients, hence the ellipticity is lost. Secondly, the localized resonance, i.e., the field explodes in some regions and remains bounded in some others as the loss goes to 0, might appear, see  \cite[Figure 2]{LaiChenZhangChanComplementary}. 
 
\medskip
In this paper, we give a proof of cloaking using complementary media  for a class of  schemes inspired by the work of Lai et al. in  \cite{LaiChenZhangChanComplementary} in the quasistatic regime. To handle the localized
resonance, we introduce the technique of removing localized singularity and apply a three spheres inequality.  The proof also uses the reflecting technique in \cite{Ng-Complementary}. To our knowledge, this work presents the first proof on cloaking using complementary media.

\begin{figure}[h!]
\begin{center}
\includegraphics[width=8cm]{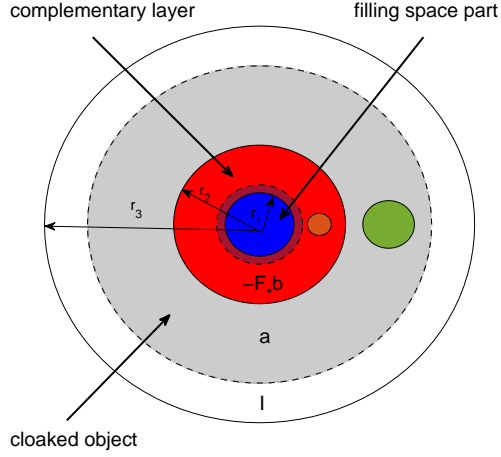}
\caption{The cloaked region in $B_{r_3} \setminus B_{r_2}$ consists of two parts: $a$  in $B_{2r_2} \setminus B_{r_2}$ (grey and green regions), which is  the original object,  and $I$ in $B_{r_3} \setminus B_{2r_2}$. The cloaking device contains two parts. The first part $-F_*b$  (the red, orange, and pink regions) in $B_{r_2} \setminus B_{r_1}$ is the complementary medium of $b$: the red part is complementary to the grey part, the orange part is the complementary to the green part, and the pink part is the complementary to $I$-layer. 
The second part $(r_3^2/r_2^2)^{d-2} I $ (the blue region)  is to fill the space  which disappears by the cancelation.} \label{fig1}
\end{center}
\end{figure}

\medskip
Let us describe how to cloak the region $B_{2r_2} \setminus B_{r_2}$ for some $r_{2}> 0$ in which the medium is characterized by a matrix $a$ using complementary media. Here and in what follows given $r>0$, $B_{r}$ denotes the ball in $\mR^d$ ($d=2$ or 3) centered at the origin of radius $r$. The assumption on the cloaked region by all means imposes no restriction  since any bounded set is a subset of such a region provided that the radius and the origin are appropriately chosen.  The idea suggested by Lai et al. in \cite{LaiChenZhangChanComplementary} (for two dimensions) is to construct a complementary medium of $a$ in $B_{r_2} \setminus B_{r_1}$ for some $0 < r_1< r_2$. Inspired by their idea, we construct a cloak  in two and three dimensions as follows. Our cloak consists of  two parts. The first one, in $B_{r_2} \setminus B_{r_1}$,  makes use of complementary media to cancel the effect of the cloaked region. The second one, in $B_{r_1}$,  is to fill the space which ``disappears" from the cancellation by the homogeneous media. 
For the first part, we slightly change the strategy in \cite{LaiChenZhangChanComplementary}. Instead of $B_{2r_2} \setminus B_{r_2}$, we consider $B_{r_3} \setminus B_{r_2}$ for some $r_3 > 0$ as the cloaked region in which the medium is given by the matrix 
\begin{equation*}
b = \left\{ \begin{array}{cl} a & \mbox{ in } B_{2 r_2} \setminus B_{r_2}, \\[6pt]
I & \mbox{ in } B_{r_3} \setminus B_{2 r_2}. 
\end{array} \right. 
\end{equation*} 
The complementary medium in $B_{r_2} \setminus B_{r_1}$ is given by 
\begin{equation*}
- \big(F^{-1}\big)_*b, 
\end{equation*}
where $F: B_{r_2} \setminus \bar B_{r_1} \to B_{r_3} \setminus \bar B_{r_2}$ is the Kelvin transform with respect to  $\partial B_{r_2}$, i.e., 
\begin{equation}\label{def-F}
F(x) = \frac{r_2^2}{|x|^2} x. 
\end{equation}
Here and in what follows we use the standard notation 
\begin{equation*}
T_*b(y) = \frac{D T  (x)  b(x) D T ^{T}(x)}{J(x)} \mbox{ where }  x =T^{-1}(y) \mbox{ and } J(x) = |\det D T(x)|, 
\end{equation*}
for a diffeomorphism $T$. It follows that 
\begin{equation}\label{cond-r}
 r_1 = r_2^2/ r_3.   
\end{equation}
Concerning the second part, the medium in $B_{r_1}$ is given by 
\begin{equation}\label{choice-A}
\Big(r_3^2/r_2^2 \Big)^{d-2} I. 
\end{equation}
The reason for this choice is condition \eqref{GF} (mentioned later in the introduction). Note that in the two dimensional case, the medium in $B_{r_1}$ is $I$, as used in \cite{LaiChenZhangChanComplementary}, while it is not $I$ in the three dimensional case.  The cloaking scheme discussed here can be extended for a large class of reflections  considered in \cite{Ng-Complementary}. The cloaking setting is illustrated in Figure~\ref{fig1}.  

To study the problem correctly, one should allow some loss in the medium and study the limit as the loss goes to 0. With the loss, the medium is characterized by $s_\delta A$ where 
\begin{equation}\label{defA}
A = \left\{ \begin{array}{cl} 
b & \mbox{ in } B_{r_3} \setminus B_{r_2}, \\[6pt]
F^{-1}_*b & \mbox{ in } B_{r_2} \setminus B_{r_1}, \\[6pt]
\Big(r_3^2/r_2^2 \Big)^{d-2} I & \mbox{ in } B_{r_1},\\[6pt]
I & \mbox{ otherwise}, 
\end{array} \right. 
\end{equation} 
and 
\begin{equation}\label{def-ss}
\epss_\delta = \left\{ \begin{array}{cl} -1 + i \delta  & \mbox{ in } B_{r_2} \setminus B_{r_1}, \\[6pt]
1 & \mbox{ otherwise}. 
\end{array} \right. 
\end{equation}
Physically, the imaginary part of $\epss_\delta A$ is the loss of the medium (more precisely the loss of  the medium of negative index in 
$B_{r_2} \setminus B_{r_1}$). Here and in what follows, we assume that,
\begin{equation}\label{cond-a}
\frac{1}{\Lambda} |\xi|^{2} \le \langle a (x) \xi, \xi \rangle \le \Lambda |\xi|^{2} \quad \forall \, \xi \in \mR^{d}, \mbox{ for a.e. } x \in B_{2r_{2}} \setminus B_{r_2}, 
\end{equation}
for some $\Lambda \ge 1$. Here and in what follows $\langle \cdot, \cdot \rangle$ denotes the Euclidean scalar product in $\mR^d$.  We require in addition that 
\begin{equation*}
b \mbox{ is Lipschitz}.
\end{equation*}

This paper deals with the bounded setting equipped  the zero Dirichlet boundary condition in the quasistatic regime. Let $\Omega$ be a smooth open subset of $\mR^d$ ($d=2, \, 3$) such that $B_{r_3} \subset \subset \Omega$. Given $f \in L^{2}(\Omega)$, let  $u_\delta, \, u \in H^1_{0}(\Omega)$ be respectively the unique solution to 
\begin{equation}\label{eq-uu-delta}
\dive (\epss_\delta A \nabla u_\delta) = f \mbox{ in } \Omega, 
\end{equation}
and 
\begin{equation}\label{eq-uu}
\Delta u = f \mbox{ in } \Omega. 
\end{equation}

Here is the main result of this paper. 

\begin{theorem}\label{thm1} Let  $d=2, \, 3$, $f \in L^{2}(\Omega)$ with $\supp f \subset \Omega \setminus  B_{r_{3}}$ and let  $u, u_\delta \in  H^1_{0}(\Omega)$ be the unique solution to \eqref{eq-uu-delta} and \eqref{eq-uu} respectively. There exists $\ell > 0$, depending only on $r_2$, $\Lambda$,  and the Lipschitz constant of $b$  such that if $r_3 > \ell r_2$ then   
\begin{equation}\label{key-point}
u_{\delta} \to u \mbox{ weakly in } H^{1} (\Omega \setminus B_{r_3}) \mbox{ as } \delta \to 0. 
\end{equation}
\end{theorem}

For an observer outside $B_{r_3}$, the medium in $B_{r_3}$ looks like the homogeneous one by \eqref{key-point} (see also \eqref{eq-uu}): one has cloaking. 

\medskip
One can verify that  medium $s_0 A$ is of reflecting complementary property, a concept introduced in \cite[Definition 1]{Ng-Complementary}, by considering  diffeomorphism  $G: \mR^{d} \setminus \bar B_{r_{3}} \to B_{r_{3}} \setminus \{0\}$ which is the Kelvin transform with respect to  $\partial B_{r_{3}}$,  i.e., 
\begin{equation}\label{def-G}
G(x) = r_{3}^{2} x/ |x|^{2}. 
\end{equation}
It is important to note that  
\begin{equation}\label{GF}
G_{*} F_{*} A = I \mbox{ in } B_{r_{3}}
\end{equation}
since $G\circ F (x) = (r_3^2/ r_2^2)x$. This is the reason behind the choice of  $A$  in  \eqref{choice-A}.  
%

\medskip
The proof of Theorem~\ref{thm1} is based on two  crucial ingredients. The first one is a  three spheres inequality. The second one is the technique of removing localized singularity  introduced in this paper \footnote{The terminology ``oscillation" has been frequently used  in the physics literature to mention the localized resonance. Here we use the standard terminology from mathematics community: ``singularity". Both of them would be fine if one looks at the  removed term  $\hat u_\delta$ defined in \eqref{def-hu} in the proof of Theorem~\ref{thm1}.}.  The proof of Theorem~\ref{thm1} is inspired by the approach in \cite{Ng-Complementary} where the reflecting technique is used. The removing singularity technique  is used in \cite{Ng-superlensing, Ng-CALR} to study superlensing and cloaking via anomalous localized resonance. 

\medskip
NIMs  have been studied extensively recently see  \cite{AmmariCiraoloKangLee, AmmariCiraoloKangLeeMilton2, AmmariCiraoloKangLeeMilton1, BouchitteSchweizer10,  KohnLu, NicoroviciMcPhedranMiltonPodolskiy1,  MiltonNicorovici,  Ng-Complementary, MinhLoc} and references therein. 
In \cite{AmmariCiraoloKangLee} and \cite{KohnLu}, Ammari et al. and Kohn et al. studied the blow up of the power ($\sim \delta \| u_\delta \|_{H^1}^2$) for a general core-shell structure with piecewise constant coefficients.  Without the core, it is shown in \cite{MinhLoc} that there is no connection between the blow up of the power and the localized resonance in general. In the works mentioned, the localized resonance was considered only in the cases where the settings are simple enough so that the separation of variables can be applied.  


\medskip
The paper is organized as follows. In Section~\ref{sect-pre}, we state and prove several useful lemmas which are used in the proof of Theorem~\ref{thm1}. The proof of Thereom~\ref{thm1} is presented in Section~\ref{sect-proof}.

\section{Preliminaries}\label{sect-pre}

In this section, we present some lemmas which will be used in the proof of Theorem~\ref{thm1}. 
The first one is on a three spheres inequality for Lipschitz matrix-valued functions. 

\begin{lemma} [Three spheres inequality] \label{lem1.1}  Let $d = 2, 3$,  $0< R_{1} < R_{2} < R_{3}$ and let $M$ be a Lipschitz matrix-valued function defined in $B_{R_{3}}$ such that $M$ is uniformly elliptic in $B_{R_3}$, 
$M(0) = I$,  and $M(x) = I$ for $x \in B_{R_3}  \setminus B_{R_2/2}$. 
Let $v \in H^{1}(B_{R_{3}})$ be a solution to  the equation
\begin{equation*}
\dive (M \nabla v ) = 0 \mbox{ in } B_{R_{3}}. 
\end{equation*}
Then for all $0 < \alpha < 1$, there exists $R_*$ depending only on $R_1$, $R_2$,  and  the Lipschitz and ellipticity constants of $M$ such that,  for $R_3 > R_*$, 
\begin{equation}\label{rem}
\|v \|_{L^{2}(B_{R_{2}})} \le C \| v\|_{L^{2}(B_{R_{1}})}^{ \alpha }  \| v\|_{L^{2}(B_{R_{3}})}^{1 - \alpha }, 
\end{equation}
for some positive constant  $C$  independent of $v$. 
\end{lemma}

\noindent{\bf Proof.} By \cite[Theorem 2.3 and (2.10)]{AlessandriniRondi}, there exists a constant $0 < \beta < 1$ which depends only on $R_1$ and $R_2$ such that  
\begin{equation}\label{part1}
\|v \|_{L^{2}(B_{R_{2}})} \le C \| v\|_{L^{2}(B_{R_{1}})}^{ \beta }  \| v\|_{L^{2}(B_{2R_{2}})}^{1 - \beta}. 
\end{equation}
Here and in what follows in this proof, $C$ denotes a positive constant which is  independent of $v$ and can change from one place to another. 
We claim that 
\begin{equation}\label{part2}
N(v, 2 R_2) \le C N(v, R_2/2)^\gamma N(v, R_3/2)^{1 - \gamma}, 
\end{equation}
where 
\begin{equation*}
N(v, r): = \|v \|_{H^{1/2}(\partial B_{r})} +\|\partial_r v \|_{H^{-1/2}(\partial B_{r})} 
\end{equation*}
and 
\begin{equation*}
\gamma = \ln \big[R_3/(4R_2) \big] \big/ \ln \big[R_3 / R_2 \big]. 
\end{equation*}
The proof of \eqref{part1} can be proceeded as follows. We only consider the case $d=2$, the case $d=3$ follows similarly. Since $\Delta v = 0$ in $B_{R_3} \setminus B_{R_2/2}$, one can represent $v$ in $B_{R_3} \setminus B_{R_2/2}$ as follows
\begin{equation*}
v(x) = \sum_{n \ge 0}  \sum_{\pm} (a_{n, \pm} r^{n, \pm} + b_n r^{-n}) e^{\pm i n \theta}
\end{equation*}
for  $a_{n, \pm}, b_{n, \pm} \in \mC$ ($n \ge 0$) with the convention $a_{0, +} = a_{0, -} = b_{0, +} = b_{0, -}$. For $R_2/2 \le r \le  R_3$, we have
\begin{equation*}
N(v, r)^2 \sim \sum_{n \ge 0}  \sum_{\pm} (n +1) \big( |a_{n, \pm}|^2 r^{2n}+ |b_{n, \pm}|^2 r^{-2n} \big). 
\end{equation*}
Here for two nonnegative quantities $\tau_1$ and $\tau_2$,  $\tau_1 \sim \tau_2$ means that $\tau_1 \le C \tau_2$ and $\tau_2 \le C \tau_1$. 
It follows that,  for $R_2/2 \le r_1 < r_2 < r_3 \le R_3$, 
\begin{equation*}
N(v, r_2) \le C N(v, r_1)^\lambda N(v, r_3)^{1 - \lambda}, 
\end{equation*}
where $\lambda = \lambda(r_1, r_2, r_3) := \ln (r_3/r_1) / \ln (r_3/ r_2)$.  We  obtain claim \eqref{part2}.

A combination of \eqref{part1} and \eqref{part2} yields 
\begin{equation}\label{part3}
\|v \|_{L^{2}(B_{R_{2}})} \le C \| v\|_{L^{2}(B_{R_{1}})}^{ \frac{\beta}{1 - \gamma (1 - \beta)} }  \| v\|_{L^{2}(B_{R_{3}})}^{\frac{(1 - \beta)(1 - \gamma)}{1 - \gamma(1 - \beta)}}. 
\end{equation}
Here we use the fact that 
\begin{equation*}
 \|v \|_{L^2(B_{2R_2})} \le C N(v, 2 R_2), \quad N(v, R_2/2) \le C \|v \|_{L^2(B_{R_2})}, \quad \mbox{ and } \quad N(v, R_3/2) \le C \|v \|_{L^2(B_{R_3})}. 
\end{equation*}
By taking $R_3$ large enough so that $\gamma$ is close to 1,   $\frac{(1 - \beta)(1 - \gamma)}{1 - \gamma(1 - \beta)}$ is close to 0. The conclusion follows. \proofend

\medskip 
The second lemma of this section is standard and its proof is left to the reader. 

\begin{lemma}\label{lem3} Let $d \ge 2$, $\alpha, \beta > 0$,  $D \subset \subset \Omega$ be two smooth open subsets of $\mR^{d}$,  and $f \in L^2(\Omega)$, $g \in H^{1/2}(\partial \Omega)$, and $h \in H^{-1/2}(\partial D)$. Assume that  $v \in H^{1}(\Omega \setminus \partial D)$ is such that 
\begin{equation}\label{v1}
\Delta v = f \mbox{ in } \Omega \setminus \partial D, \quad v = 0 \mbox{ on } \partial \Omega, 
\end{equation}
and
\begin{equation}\label{v2}
\quad [v] = g \mbox{ on } \partial D, \quad \mbox{ and } \quad [\partial_r v] = h \mbox{ on } \partial D.  
\end{equation}
We have
\begin{equation}\label{v4}
\|v \|_{H^{1}(\Omega \setminus \overline D)} + \|v \|_{H^{1}(D)} \le C  \Big( \| f\|_{L^2(\Omega)} + \| g \|_{H^{1/2}(\partial D)}  +  \| h \|_{H^{-1/2}(\partial D)}  \Big), 
\end{equation}
for some positive constant $C$ independent of $f$, $g$,  and $h$. 
\end{lemma}

Here and in what follows, $v\big|_{+}$, $v\big|_{-}$, and $[v]$  on $\partial D$ denote the trace of a function $v$ from the exterior, interior of $D$, and $v\big|_{+} - v\big|_{-}$ on $\partial D$ respectively. Similarly, the notations  ${\cal M} \nabla v \cdot \eta\big|_{+}$, ${\cal M} \nabla v \cdot \eta \big|_{-}$, and $[{\cal M} \nabla v \cdot \eta]$  are used on $\partial D$ for an appropriate matrix ${\cal M}$ and an appropriate function $v$, where $\eta$ is the normal unit vector directed to the exterior of $D$.

\medskip 
The following lemma which is a consequence of  \cite[Lemma 4]{Ng-Complementary} is  used  in the proof of Theorem~\ref{thm1}. 

\begin{lemma}\label{lem-TO} Let $d \ge 2$,  $0 < R_1 < R_2 < R_3$ with $R_3 = R_2^2/ R_1$, $a \in [L^\infty(B_{R_3 \setminus R_2})]^{d \times d}$ be a matrix valued function, and $K:B_{R_2} \setminus \bar B_{R_1} \to B_{R_3} \setminus \bar B_{R_2}$ be the Kelvin transform with respect to $\partial B_{R_2}$, i.e., 
\begin{equation*}
K(x) = R_2^2 x/ |x|^2. 
\end{equation*}
For $v \in H^1(B_{R_2} \setminus B_{R_1})$, define $w = v \circ F^{-1}$. Then 
\begin{equation*}
\dive (a \nabla v) = 0 \mbox{ in }  B_{R_2} \setminus B_{R_1}
\end{equation*}
if and only if
\begin{equation*}
\dive (K_*a \nabla w)  = 0 \mbox{ in } B_{R_3} \setminus B_{R_2}.  
\end{equation*}
Moreover, 
\begin{equation*}
w = v \quad \mbox{ and } \quad K_*a \nabla w \cdot \eta = - a \nabla v \cdot \eta  \mbox{ on } \partial B_{R_2}.    
\end{equation*}
\end{lemma}

\section{Proof of Theorem~\ref{thm1}} \label{sect-proof}

Let ${\cal A}$ be a matrix valued function defined in $B_{r_3}$ such that ${\cal A} = A = b$ in $B_{r_3} \setminus B_{r_2}$, ${\cal A}(0) = I$,    and ${\cal A}$ is Lipschitz. Take $\ell  > 4$ large enough  such that \eqref{rem} holds for 
\begin{equation}\label{alpha}
\alpha =  2/3, 
\end{equation}
with  $M= {\cal A}$, $R_1 = r_2$, $R_2 = 4r_2$, and $R_3 = r_3$. Then 
\begin{equation}\label{main-point}
\|v \|_{L^2(B_{4 r_2})} \le \|v \|_{L^2(B_{r_2})}^{\alpha} \|v \|_{L^2(B_{r_3})}^{1 - \alpha},   
\end{equation}
for $v \in H^1(B_{r_3})$ which satisfies $\dive ( {\cal A} \nabla v ) = 0$ in $B_{r_3}$. 

\medskip
Multiplying \eqref{eq-uu-delta} by $\bar u_{\delta}$ and integrating on $\Omega$, we obtain 
\begin{equation*}
\int_{\Omega} \langle s_\delta A \nabla u_{\delta}, \nabla u_{\delta} \rangle = - \int_{\Omega} f  \bar u_{\delta}. 
\end{equation*}
Considering the imaginary part and the real part, we have
\begin{equation}\label{var1}
\int_{\Omega} |\nabla u_{\delta}|^{2} \le C\delta^{-1} \| f\|_{L^2(\Omega)} \|u_{\delta} \|_{L^2(\Omega \setminus B_{r_3})}.
 \end{equation}
In this proof, $C$ denotes a positive constant changing from one place to another but independent of $\delta$ and $f$.  Since $u_{\delta} = 0$ on $\partial \Omega$, it follows that 
\begin{equation*}
\| u_{\delta}\|_{L^{2}(\Omega)} \le C \| \nabla u_{\delta} \|_{L^{2}(\Omega)}. 
\end{equation*}
We derive from \eqref{var1} that 
\begin{equation}\label{id1}
\| u_{\delta} \|_{H^{1}(\Omega)} \le C \delta^{-1/2} \| f\|_{L^2(\Omega)}^{1/2} \| u_\delta \|_{L^2(\Omega \setminus B_{r_3})}^{1/2}
\end{equation}
and
\begin{equation}\label{id1-1}
\| u_{\delta} \|_{H^{1}(\Omega)} \le C \delta^{-1}\| f\|_{L^2(\Omega)}.
\end{equation}
As in \cite{Ng-Complementary}, let $u_{1, \delta}$ be the reflection of $u_{\delta}$ through $\partial B_{r_{2}}$ by $F$,  i.e., 
\begin{equation}\label{def-u1}
u_{1, \delta } = u_{\delta} \circ F^{-1} \mbox{ in } \mR^{d} \setminus \bar B_{r_{2}}. 
\end{equation}
Applying Lemma~\ref{lem-TO} and using the fact that $F_*A = A$ in $B_{r_3} \setminus B_{r_2}$, we obtain 
\begin{equation}\label{TO-1}
u_{1, \delta} = u_{\delta} \big|_{+}   \mbox{ on } \partial B_{r_{2}} \quad \mbox{ and } \quad  (1 - i \delta ) A \nabla u_{1, \delta} \cdot \eta = A \nabla u_{ \delta} \cdot \eta  \big|_{+} \mbox{ on } \partial B_{r_{2}}. 
\end{equation}
Define  $\bU_{\delta}$ in $B_{r_{3}}$ as follows
\begin{equation}\label{def-bU}
\bU_{\delta} = \left\{\begin{array}{cl} u_{\delta} - u_{1, \delta} & \mbox{ in } B_{r_{3}} \setminus B_{r_{2}}, \\[6pt]
0 & \mbox{ in }  B_{r_{2}}.
\end{array} \right.
\end{equation}
Using the fact that $F_*A = A$ in $B_{r_3} \setminus B_{r_2}$ and applying Lemma~\ref{lem-TO}, we derive from \eqref{def-u1}, \eqref{TO-1}, and \eqref{def-bU} that  $\bU_{\delta} \in H^{1}(B_{r_{3}})$, 
\begin{equation}\label{eq-bU}
\dive ({\cal A} \nabla \bU_{\delta} )  = 0 \mbox{ in } B_{r_{3}} \setminus \partial B_{r_{2}}, \quad \mbox{ and } \quad  [{\cal A} \nabla \bU_{\delta} \cdot \eta] = - \frac{i \delta}{1 - i \delta } {\cal A} \nabla u_{\delta} \cdot \eta  \big|_{+}\mbox{ on } \partial B_{r_{2}}.   
\end{equation}
 It is clear from \eqref{id1}, \eqref{def-u1},  and \eqref{def-bU} that 
\begin{equation}\label{est-bU}
\| \bU_\delta \|_{H^{1}(B_{r_3})} \le C \delta^{-1/2} \| f\|_{L^2(\Omega)}^{1/2} \| u_\delta \|_{L^2(\Omega \setminus B_{r_3})}^{1/2}. 
\end{equation}
Let $w_{\delta} \in H^{1}_{0}(B_{r_{3}})$ be the unique solution to 
\begin{equation}\label{def-w}
\dive ({\cal A} \nabla w_{\delta}) = 0 \mbox{ in } B_{r_{3}} \setminus \partial B_{r_{2}} \quad \mbox{ and } \quad [{\cal A} \nabla w_{\delta}  \cdot \eta] = - \frac{i \delta}{1 - i \delta } {\cal A} \nabla u_{\delta} \cdot \eta \big|_{+}\mbox{ on } \partial B_{r_{2}}. 
\end{equation}
Then 
\begin{equation}\label{id2}
\|w_{\delta} \|_{H^{1}(B_{r_{3}})} \le C \delta \| A \nabla u_{\delta } \cdot \eta \big|_{+}\|_{H^{-1/2}(\partial B_{r_{2}})}. 
\end{equation}
Since 
\begin{equation*}
\| A \nabla u_{\delta } \cdot \eta \big|_{+}\|_{H^{-1/2}(\partial B_{r_{2}})} \le  C \| u_\delta \|_{H^1(\Omega)}, 
\end{equation*}
it follows from  \eqref{id1} and \eqref{id2} that 
\begin{equation}\label{id3}
\| w_{\delta}\|_{H^{1}(B_{r_{3}})} \le C \delta^{1/2} \| f\|_{L^2(\Omega)}^{1/2} \| u_\delta \|_{L^2(\Omega \setminus B_{r_3})}^{1/2}. 
\end{equation}
Define 
\begin{equation}\label{def-bV}
\bV_{\delta} = \bU_{\delta} - w_{\delta} \mbox{ in } B_{r_{3}}.  
\end{equation}
Then $\bV_{\delta} \in H^{1}(B_{r_{3}})$  is a solution to 
\begin{equation}\label{eq-V}
\dive ({\cal A} \nabla \bV_{\delta}) = 0 \mbox{ in } B_{r_{3}}. 
\end{equation}
Using \eqref{main-point}, we obtain 
\begin{equation}\label{id4}
\|\bV_{\delta} \|_{L^{2}(B_{4 r_2})} \le C \|\bV_{\delta} \|_{L^{2}(B_{r_{2}})}^{\alpha} \|\bV_{\delta} \|_{L^{2}(B_{r_{3}})}^{1 - \alpha}.
\end{equation}
From \eqref{def-bU} and \eqref{def-bV}, we have
\begin{equation*}
\bV_{\delta} = - w_{\delta} \mbox{ in } B_{r_{2}}. 
\end{equation*}
We derive from  \eqref{id3}  that
\begin{equation}\label{id5}
\| \bV_{\delta}\|_{L^{2}(B_{r_{2}})} \le C \delta^{1/2} \| f\|_{L^2(\Omega)}^{1/2} \| u_\delta \|_{L^2(\Omega \setminus B_{r_3})}^{1/2}. 
\end{equation}
On the other hand, from \eqref{est-bU}, \eqref{id3}, and \eqref{def-bV}, we have
\begin{equation}\label{id6}
\| \bV_{\delta}\|_{L^{2}(B_{r_{3}})} \le C \delta^{-1/2} \| f\|_{L^2(\Omega)}^{1/2} \| u_\delta \|_{L^2(\Omega \setminus B_{r_3})}^{1/2}. 
\end{equation}
A combination of  \eqref{id4}, \eqref{id5}, and \eqref{id6} yields
\begin{equation*}
\|\bV_{\delta} \|_{L^{2}(B_{4 r_2})}^{2}  \le C \delta^{(2 \alpha - 1)} \|f\|_{L^2(\Omega)} \| u_\delta \|_{L^2(\Omega \setminus B_{r_3})}.  
\end{equation*}
Since $\Delta \bV_{\delta} = 0 $ in $B_{4 r_{2} } \setminus B_{2 r_{2}}$, it follows that 
\begin{equation}\label{id6-1}
\|\bV_{\delta} \|_{H^{1/2}(\partial B_{3 r_2})}^{2}  +  \|\partial_r \bV_{\delta} \|_{H^{-1/2}(\partial B_{3 r_2})}^{2}  \le C \delta^{(2 \alpha - 1)} \|f\|_{L^2(\Omega)} \| u_\delta \|_{L^2(\Omega \setminus B_{r_3})}. 
\end{equation}
A combination of \eqref{def-bU}, \eqref{id3}, \eqref{def-bV}, and \eqref{id6-1} implies
\begin{equation}\label{id7}
\| u_{\delta} - u_{1, \delta} \|_{H^{1/2}(\partial B_{3 r_2}) }^{2}  + \| \partial_{r} u_{\delta} - \partial_{r} u_{1, \delta} \|_{H^{-1/2}(\partial B_{3 r_2}) }^{2}
\le C \delta^{2 \beta}  \|f\|_{L^2(\Omega)} \| u_\delta \|_{L^2(\Omega \setminus B_{r_3})}, 
\end{equation}
where, by  \eqref{alpha},  
\begin{equation}
0<  \beta : = (2 \alpha - 1)/2  = 1/ 6 < 1/ 2. 
\end{equation}
As in \cite{Ng-Complementary}, let $u_{2, \delta}$ be the reflection of $u_{1, \delta}$ through $\partial B_{r_{3}}$ by $G$, i.e., 
\begin{equation*}
u_{2, \delta} = u_{1, \delta} \circ G^{-1}  \mbox{ in } B_{r_{3}}. 
\end{equation*}
We have
\begin{equation}\label{u1u2Delta}
\Delta u_{1, \delta} = 0 \mbox{ in } B_{r_{3}} \setminus B_{3 r_2} \quad \mbox{ and }  \quad \Delta u_{2, \delta} = 0 \mbox{ in } B_{r_{3}}, 
\end{equation}
by \eqref{GF}. 

\medskip 
We next consider the case $d=2$ and the case $d=3$ separately. 

\medskip
\noindent \underline{Case 1:} $d=2$. From \eqref{u1u2Delta}, one can represent $u_{1, \delta}$ and $u_{2, \delta}$ as follows 
\begin{equation}\label{eq-u1}
u_{1, \delta} =  c_{0} + d_{0} \ln r +  \sum_{n = 1}^\infty \sum_{\pm} (c_{n, \pm} r^{n} + d_{n, \pm} r^{-n} ) e^{ \pm i n \theta } 
\mbox{ in } B_{r_{3}} \setminus B_{3 r_2}, 
\end{equation}
and 
\begin{equation}\label{eq-u2}
u_{2, \delta} = e_{0} + \sum_{n = 1}^\infty \sum_{\pm} e_{n, \pm} r^{n}  e^{\pm i n \theta}  \mbox{ in } B_{r_{3}},  
\end{equation}
for $c_0, d_0, e_0, c_{n, \pm}, d_{n, \pm}, e_{n, \pm} \in \mC$ ($n \ge 1$). 
We derive from  \eqref{id1} that 
\begin{equation}\label{est-e}
|e_0|^2 + \sum_{n = 1}^\infty \sum_{\pm} n |e_{n, \pm}|^{2} r_{3}^{2n} \le C \delta^{-1} \| f\|_{L^2(\Omega)}  \| u_\delta \|_{L^2(\Omega \setminus B_{r_3})}.
\end{equation}
Using the fact that $G_*F_*A = I$ and applying Lemma~\ref{lem-TO},  we have
\begin{equation}\label{transmission-r3}
u_{1, \delta} = u_{2, \delta} \mbox{ on } \partial B_{r_{3}} \quad  \mbox{ and }  \quad (1 - i \delta )\partial_{\eta} u_{1, \delta} \big|_{-}= \partial_{\eta} u_{2, \delta} \mbox{ on } \partial B_{r_{3}}.
\end{equation}
A combination of \eqref{eq-u1}, \eqref{eq-u2}, and \eqref{transmission-r3} yields
\begin{equation}\label{eq-abe*}
c_{n, \pm} r_3^{n} + d_{n, \pm} r_3^{-n} = e_{n, \pm} r_3^n,  \quad c_{n, \pm} r_3^{n} - d_{n, \pm} r_3^{-n} = \frac{1}{1 - i \delta} e_{n, \pm} r_3^n  \mbox{ for } n \ge 1, 
\end{equation}
and 
\begin{equation}\label{eq-abe-1}
c_{0} = e_{0}, \quad d_{0} = 0. 
\end{equation}
We derive from \eqref{eq-abe*} that 
\begin{equation}\label{eq-abe}
c_{n, \pm} = \frac{2 - i \delta}{2(1 - i \delta)} e_{n, \pm} \quad  \mbox{ and } \quad d_{n, \pm} = \frac{- i \delta}{2 (1  - i \delta) }  r_{3}^{2n}  e_{n, \pm} \quad \mbox{ for } n \ge 1, 
\end{equation}
From \eqref{eq-abe-1} and \eqref{eq-abe}, we obtain
\begin{equation}\label{eq-cde}
u_{1, \delta} - u_{2, \delta} = \sum_{n =1 }^\infty \sum_{\pm} \frac{i \delta}{2 (1 - i \delta)} \Big( r^n - \frac{r_3^{2n}}{r^n} \Big) e_{n, \pm} e^{\pm i n \theta}. 
\end{equation}

We now introduce the technique of removing localized singularity. 
Set 
\begin{equation}\label{def-hu}
\hu_{\delta} : = \sum_{n = 1}^\infty \sum_{\pm}  \frac{- i \delta}{2 (1  - i \delta) } e_{n, \pm} \frac{r_{3}^{2n}}{r^{n}} e^{\pm i n \theta} \mbox{ in } \mR^{2} \setminus B_{3 r_2}. 
\end{equation}
Define \footnote{We remove $\hat u_\delta$ from $u_\delta$ in $B_{r_3} \setminus B_{3 r_2}$. Function $\hat u_\delta$ contains  high modes and creates a trouble for estimating  $u_\delta - u_{2, \delta}$ on $\partial B_{3 r_2}$ (to obtain an estimate for $u_\delta$). However this term can be negligible  for large $|x| = r_3$ since $r^{-n}$ is small for large $r$.  This is the reason for the choice of $\hat u_\delta$ and this also explains the terminology: ``removing localized singularity".} 
\begin{equation*}
W_{\delta} = \left\{\begin{array}{cl} u_\delta & \mbox{ in } \Omega \setminus B_{r_3}, \\[6pt]
u_{\delta} -  \hu_{\delta} & \mbox{ in } B_{r_3} \setminus B_{3 r_2},  \\[6pt]
u_{2, \delta} & \mbox{ in } B_{3 r_2}. 
\end{array} \right.
\end{equation*}
It is clear from \eqref{u1u2Delta} that 
\begin{equation}\label{eq-W}
\Delta W_{\delta} = f \mbox{ in } \Omega \setminus (\partial B_{r_3} \cup \partial B_{3 r_2}). 
\end{equation}

We claim that 
\begin{equation}\label{claim-hu}
\| \hu_{\delta} \|_{H^{1/2}(\partial B_{r_3})} + \|\partial_r \hu_{\delta} \|_{H^{-1/2}(\partial B_{r_3})} = o(1) \| f\|_{L^2(\Omega)} + o(1) \| u_\delta\|_{L^2(\Omega \setminus B_{r_3})}. 
\end{equation}
Here and in what follows we use the standard notation: $o(1)$ denotes a quantity which converges to $0$ as $\delta \to 0$. Indeed, from  \eqref{est-e}, we have
\begin{equation}\label{DW3}
\sum_{n = 1}^\infty \sum_{\pm} n \delta^2  |e_{n, \pm}|^2 r_3^{2n} \le  C  \delta  \| f\|_{L^2(\Omega)}  \| u_\delta \|_{L^2(\Omega \setminus B_{r_3})} \le o(1) \| f\|_{L^2}^2 + o(1)\| u_\delta \|_{L^2(\Omega \setminus B_{r_3})}^2. 
\end{equation}
Hence  \eqref{claim-hu} follows. We derive from \eqref{claim-hu} that 
\begin{equation}\label{claim-hu-1}
\| [W_{\delta}] \|_{H^{1/2}(\partial B_{r_3})} + \|[\partial_r W_{\delta}] \|_{H^{-1/2}(\partial B_{r_3})} = o(1) \| f\|_{L^2(\Omega)} + o(1) \| u_\delta\|_{L^2(\Omega \setminus B_{r_3})}. 
\end{equation}

\medskip
We next consider the transmission conditions for $W_\delta$ on $\partial B_{3 r_2}$. Since 
\begin{equation*}
[W_\delta] =  u_\delta - \hu_\delta - u_{2, \delta} \quad \mbox{ on } \partial B_{3 r_2}, 
\end{equation*}
it follows  from  \eqref{eq-cde} and  \eqref{def-hu} that
\begin{equation*}
[W_\delta] = (u_\delta - u_{1, \delta}) +  \sum_{n = 1}^\infty \sum_{\pm} \frac{i \delta}{2 (1 - i \delta)}  (3 r_2)^n e_{n, \pm} e^{i n \theta}   \quad \mbox{ on } \partial B_{3 r_2}. 
\end{equation*}
This implies 
\begin{equation}\label{T-1}
C \| [W_{\delta}] \|_{H^{1/2}(\partial B_{3 r_2})}^2  \le \| u_\delta - u_{1, \delta} \|_{H^{1/2}(\partial B_{3 r_2})}^2 \\[6pt] + \sum_{n = 1}^\infty \sum_{\pm} \delta^{2} n |e_{n, \pm}|^2 r_3^{2n}\frac{(3 r_2)^{2n}}{r_3^{2n}}. 
\end{equation}
We derive from \eqref{id7} and  \eqref{est-e} that 
\begin{align*}
C \| [W_{\delta}] \|_{H^{1/2}(\partial B_{3 r_2})}^2 \le &   \delta^{2 \beta} \|f\|_{L^2(\Omega)} \| u_\delta \|_{L^2(\Omega \setminus B_R)}  +  \delta \| f\|_{L^2(\Omega)} \| u_\delta \|_{L^2(\Omega \setminus B_{r_3})}; 
\end{align*}
which yields
\begin{equation}\label{T-3}
| [W_{\delta}] \|_{H^{1/2}(\partial B_{3 r_2})}^{2} 
=  o(1) \| f\|_{L^2( \Omega)}^{2} + o(1) \| u_\delta \|_{L^2(\Omega \setminus B_{r_3})}^2. 
\end{equation} 
Similarly, 
\begin{equation}\label{T-4}
 \| [\partial_{r} W_{\delta}] \|_{H^{-1/2}(\partial B_{3 r_2})}^{2} =  o(1) \| f\|_{L^2(\Omega)}^{2} + o(1) \| u_\delta \|_{L^2(\Omega \setminus B_{r_3})}^2. 
\end{equation}
Applying Lemma~\ref{lem3} and using \eqref{claim-hu-1},  \eqref{T-3}, and \eqref{T-4}, we have
\begin{equation*}
\|W_\delta \|_{H^1\big(\Omega \setminus (\partial B_{r_3} \cup \partial B_{3 r_2}) \big)} \le C \| f\|_{L^2}, 
\end{equation*}
for small $\delta$.

\medskip
Without loss of generality,  one may assume that $W_{\delta} \to W \mbox{ weakly in } H^{1}\big(\Omega \setminus (\partial B_{r_3} \cup \partial B_{3 r_2}) \big)$ as $\delta \to 0$. From \eqref{claim-hu-1}, \eqref{T-3}, and \eqref{T-4}, we have
\begin{equation*}
W \in H^1_0(\Omega) \quad \mbox{ and } \quad  \Delta W = f \mbox{ in } \Omega. 
\end{equation*}
Hence $W = u$. Since the limit $W$ is unique,  the convergence holds for the whole family $(W_{\delta})$ as $\delta \to 0$. The proof is complete in two dimensions.

\medskip
\noindent \underline{Case 2:} $d=3$.  The proof in the three dimensional case follows similarly as the one in the two dimensional  case. We just note here that, in three dimensions, $u_{1, \delta}$ and $u_{2, \delta}$ can be represented, by \eqref{u1u2Delta}, as follows
\begin{equation*}
u_{1, \delta} = c_0 + \frac{d_0}{r} +  \sum_{n=1}^\infty \sum_{k=-n}^n (c_{n, k} r^n + d_{n, k} r^{-n - 1}) Y^k_n(x/|x|) \quad \mbox{ in } B_{r_3} \setminus B_{r_2}
\end{equation*}
and 
\begin{equation*}
u_{2, \delta} = e_0 + \sum_{n=1}^\infty \sum_{k=-n}^n e_{n, k} r^n  Y^k_n(x/|x|) \quad \mbox{ in } B_{r_3}. 
\end{equation*}
 The proof is complete. \proofend

\bigskip

\noindent{\bf Acknowledgment: } The author thanks Graeme Milton for pointing out the localized resonance from simulations in \cite{LaiChenZhangChanComplementary}. He also thanks Luca Rondi for enlightening discussions on the three spheres inequality.

\providecommand{\bysame}{\leavevmode\hbox to3em{\hrulefill}\thinspace}
\providecommand{\MR}{\relax\ifhmode\unskip\space\fi MR }
\providecommand{\MRhref}[2]{%
  \href{http://www.ams.org/mathscinet-getitem?mr=#1}{#2}
}
\providecommand{\href}[2]{#2}


\end{document}